\documentclass{elsart}
\usepackage{amssymb}

\renewcommand{\bar}[1]{\overline{#1}}

\usepackage{indentfirst}
\usepackage{psfig,color}
\usepackage{epsfig}
\usepackage{epsf}
\usepackage{graphicx}

\journal{J. Differential Equations }

\begin{document}

\begin{frontmatter}
\title{Special-series solution of the first-order linear vector differential equation}

\author[pku]{Xin-Bing Huang \corauthref{cor}}
\corauth[cor]{Corresponding author.} \ead{huangxb@pku.edu.cn}
\address[pku]{Department of Physics, Peking University, Beijing 100871, China}

\begin{abstract}
A special series is introduced in this paper to yield solution of
the first-order linear vector differential equation. It is proved
that if the differential equation satisfied by the first term of
this series can be solved exactly, then other terms can be
determined by the method of variation of parameters. We point out
that the special series will be the solution of the first-order
linear vector differential equation if the infinite special series
converges. An illustrative example has been given to outline the
procedure of our method.
\end{abstract}

\end{frontmatter}

\par

\section{Introduction}

The content of linear systems constitutes a large and very
important part in the theory of ordinary differential equations.
From the early days of ordinary differential equations the subject
of linear systems has been an area of great theoretical research
and practical applications, and it continues to be so today. In
physics, many dynamical systems or material fields can be treated
as linear systems. For instance, the evolution equation of a
spin-$\frac{1}{2}$ fermion in Robertson-Walker space-time deduced
from the covariant Dirac equation is a linear system~\cite{hua05}.
In this paper we shall restrict our attention to linear systems of
$n$ differential equations in $n$ unknown functions only.

Almost in all books on ordinary differential equations, the
authors always present a symbolic operator method for solving
linear systems with constant coefficients. After introducing the
so-called differential operators, the standard procedure for
solving linear systems with constant coefficients will be
described in any such book~\cite{ros80}. Therefore, we shall only
focus on the linear systems whose coefficients are functions.

When the system of linear differential equations is of the type
that have variable coefficients, the method of solving it depends
upon the concrete form of the coefficients and the type of
differential equations. One can try to solve this kind of linear
systems by acquiring a higher order differential equation which
can be solved~\cite{wg00,ww27}, or to solve them approximately by
various methods such as numerical methods, series methods, and
graphical methods. In this paper We shall introduce a relatively
general and systematic method for solving linear systems with
variable coefficients.

The present paper is devoted to a method for obtaining solution of
the first-order linear vector differential equation in special
series form. We shall decompose the linear vector differential
equation into a group of correlated vector differential equations.
After this, we shall proceed to explain the method of special
series. We shall then discuss the convergence of infinite special
series, and point out the solution. Finally, we shall give an
illustrative example to clarify the procedure of special-series
method.

\section{Decomposing the linear vector
differential equation  }

We consider the normal form of linear system of $n$ first-order
differential equations in $n$ unknown functions $y_1$,
$y_2$,$\ldots$,$y_n$. It is well known that this system is of the
form
\begin{equation}
\label{dieq}  {\bf{y}}^{\prime}= {\bf{A}}(x){\bf{y}}+ {\bf{F}}(x)
~,
\end{equation}
where we have introduced the notation ${ {\bf{y}}}^{\prime}\equiv
\frac{d {\bf{y}}}{d x}$, the vectors ${\bf{y}}$ and ${\bf{F}}(x)$
are defined respectively by
\begin{equation}
\label{yf} {\bf{y}}=\left(
\begin{array}{c}
y_1
\\
y_2
\\
\vdots \\
 y_n
\end{array}\right)~~~~~~~~{\rm and}~~~~~~~~{\bf{F}}(x)=\left(
\begin{array}{c}
f_{1}(x)
\\
f_{2}(x)
\\
\vdots
\\
f_{n}(x)
\end{array}\right)
\end{equation}
and the matrix function ${\bf{A}}(x)$ being defined by
\begin{equation}
\label{adef} {\bf{A}}(x)=\left(
\begin{array}{cccc}
a_{11}(x) & a_{12}(x) & \cdots & a_{1n}(x)
\\
a_{21}(x) & a_{22}(x) & \cdots & a_{2n}(x)
\\
\vdots & \vdots & ~ &   \vdots
\\
a_{n1}(x) & a_{n2}(x) & \cdots & a_{nn}(x)
\end{array}\right)~.
\end{equation}
We shall assume that all of the functions defined by $f_{i}(x)$,
$i=1,2,\ldots,n$, and $a_{ij}(x)$,
$i=1,2,\ldots,n,j=1,2,\ldots,n$, are continuous on a real interval
$[a, b]$. If all $f_{i}(x)=0$, $i=1,2,\ldots,n$, for all $x$, then
the system (\ref{dieq}) is called homogeneous. Otherwise, the
system is called nonhomogeneous.

Hence the normal form of linear system (\ref{dieq}) is a
nonhomogeneous linear vector differential equation, and the
corresponding homogeneous linear vector differential equation
reads
\begin{equation}
\label{dieqhom}  {\bf{y}}^{\prime}= {\bf{A}}(x){\bf{y}} ~.
\end{equation}
Mathematician has developed the method of variation of parameters
for finding a solution of the equation (\ref{dieq}), assuming that
one has known a fundamental matrix of the corresponding equation
(\ref{dieqhom}). Therefore, in this paper, we will focus on
explaining the method of special series for finding the
fundamental matrix of the homogeneous linear vector differential
equation (\ref{dieqhom}).

We shall assume that ${\bf{A}}(x)$ is a nonsingular $n\times n$
matrix function. Then one can divide ${\bf{A}}(x)$ into two
$n\times n$ nonsingular matrices ${\bf{H}}(x)$ and ${\bf{Z}}(x)$.
It is stressed that the elements of matrix function ${\bf{Z}}(x)$
must be finite on the real interval $[a,b]$. That is,
\begin{equation}
\label{sepa} {\bf{A}}(x)= {\bf{H}}(x)+ {\bf{Z}}(x) ~.
\end{equation}
Here ${\bf{H}}(x)$ can also be a constant matrix. For brevity, we
use the common sign ${\bf{H}}(x)$ to denote both a matrix function
defined on all real $x\in [a,b]$ and a constant matrix. Therefore,
the equation (\ref{dieqhom}) becomes
\begin{equation}
\label{dieqhom2}  {\bf{y}}^{\prime}= {\bf{H}}(x){\bf{y}}+
{\bf{Z}}(x){\bf{y}} ~.
\end{equation}

To find the fundamental matrix of the above equation, we shall
assume that the vector function ${\bf{y}}$ is a special series
written as follows
\begin{equation}
\label{infsum}  {\bf{y}}=\sum^{\infty}_{i=0} {\bf{y}}_{i} ~.
\end{equation}
And we shall assume that the first term ${\bf{y}}_{0}$ of this
series satisfies
\begin{equation}
\label{y0}  {\bf{y}}^{\prime}_{0}= {\bf{H}}(x){\bf{y}}_{0}~.
\end{equation}

Under the assumption that ${\bf{y}}$ is the sum of the first two
terms, inserting ${\bf{y}}={\bf{y}}_{0}+{\bf{y}}_{1}$ into the
equation (\ref{dieqhom2}) yields
\begin{equation}
\label{y0y1}  {\bf{y}}^{\prime}_{0}+{\bf{y}}^{\prime}_{1}=
{\bf{H}}(x)({\bf{y}}_{0}+{\bf{y}}_{1})+
{\bf{Z}}(x)({\bf{y}}_{0}+{\bf{y}}_{1})~.
\end{equation}
With the help of the equation (\ref{y0}), and ignoring the term
${\bf{Z}}(x){\bf{y}}_{1}$ in the equation (\ref{y0y1}), we obtain
\begin{equation}
\label{y1}  {\bf{y}}^{\prime}_{1}= {\bf{H}}(x){\bf{y}}_{1}+
{\bf{Z}}(x){\bf{y}}_{0}~.
\end{equation}
Under the assumption that
${\bf{y}}={\bf{y}}_{0}+{\bf{y}}_{1}+{\bf{y}}_{2}$, inserting
${\bf{y}}_{0}+{\bf{y}}_{1}+{\bf{y}}_{2}$ into the equation
(\ref{dieqhom2}) yields
\begin{equation}
\label{y0y1y2}
{\bf{y}}^{\prime}_{0}+{\bf{y}}^{\prime}_{1}+{\bf{y}}^{\prime}_{2}=
{\bf{H}}(x)({\bf{y}}_{0}+{\bf{y}}_{1}+{\bf{y}}_{2})+
{\bf{Z}}(x)({\bf{y}}_{0}+{\bf{y}}_{1}+{\bf{y}}_{2})~.
\end{equation}
With the help of the equations (\ref{y0}), (\ref{y1}), and
ignoring the term ${\bf{Z}}(x){\bf{y}}_{2}$ in the above equation,
we obtain
\begin{equation}
\label{y2}  {\bf{y}}^{\prime}_{2}= {\bf{H}}(x){\bf{y}}_{2}+
{\bf{Z}}(x){\bf{y}}_{1}~.
\end{equation}
Under the assumption that ${\bf{y}}=\sum_{j=0}^{m} {\bf{y}}_{j}$,
inserting $\sum_{j=0}^{m} {\bf{y}}_{j}$ into the equation
(\ref{dieqhom2}) yields
\begin{equation}
\label{ysum}  \sum_{j=0}^{m}{\bf{y}}^{\prime}_{j}=
{\bf{H}}(x)\sum_{j=0}^{m}{\bf{y}}_{j}+
{\bf{Z}}(x)\sum_{j=0}^{m}{\bf{y}}_{j}~.
\end{equation}
Using the same method as before, and ignoring the term
${\bf{Z}}(x){\bf{y}}_{m}$ in the above equation, we acquire
\begin{equation}
\label{ym}  {\bf{y}}^{\prime}_{m}= {\bf{H}}(x){\bf{y}}_{m}+
{\bf{Z}}(x){\bf{y}}_{m-1}~.
\end{equation}

Through the foregoing logic chain, we can draw the following
conclusion: in the case of homogeneous linear vector differential
equation (\ref{dieqhom}), after expanding the vector of functions
${\bf{y}}$ into an infinite series, at the same time, dividing the
coefficient matrix into two nonsingular matrices, we can obtain a
recursion series, each term of which satisfies
\begin{equation}
\label{con1}\left\{ \begin{array}{l} {\bf{y}}^{\prime}_{j}=
{\bf{H}}(x){\bf{y}}_{j}~~~~~~~~~~~~~~~~~~~~~~~~~~j=0
\\
{\bf{y}}^{\prime}_{j}= {\bf{H}}(x){\bf{y}}_{j}+
{\bf{Z}}(x){\bf{y}}_{j-1}~~~~~~~~~~~j=1,2,3,\ldots
\end{array} \right. ~.
\end{equation}
In the case of nonhomogeneous linear vector differential equation
(\ref{dieq}), we can similarly define a special series, each term
of which satisfies
\begin{equation}
\label{con2}\left\{ \begin{array}{l} {\bf{y}}^{\prime}_{j}=
{\bf{H}}(x){\bf{y}}_{j}+{\bf{F}}(x)~~~~~~~~~~~~~~~~j=0
\\
{\bf{y}}^{\prime}_{j}= {\bf{H}}(x){\bf{y}}_{j}+
{\bf{Z}}(x){\bf{y}}_{j-1}~~~~~~~~~~~j=1,2,3,\ldots
\end{array} \right. ~.
\end{equation}
From the theoretical point of view, we can say, the above
equations demonstrate that one can obtain every terms in the
recursion series ${\bf{y}}=\sum^{\infty}_{i=0} {\bf{y}}_{i}$ if
and only if the first term can be solved explicitly.


\section{The solution of each term of special series}
The theory of linear systems has proven that if the vector
functions ${\Phi}_{1}$, ${\Phi}_{2}$, $\ldots$ ,${\Phi}_{i}$ are
$i$ solutions of (\ref{y0}) and $c_1$, $c_2$, $\ldots$, $c_i$ are
$i$ numbers, then the vector function
\begin{equation}\label{linsum}
 {\Phi}=\sum_{k=1}^{i}c_{k}{\Phi}_{k}~,
\end{equation}
is also a solution of (\ref{y0}).

We shall be concerned with $n$ vector functions, and we shall use
the following common notation for the $n$ vector functions in the
following discussion. We let ${\Phi}_{1}$, ${\Phi}_{2}$, $\ldots$
,${\Phi}_{n}$ be the $n$ vector functions defined respectively by
\begin{equation}\label{nvecdef}
{\Phi}_{1}(x)=\left(
\begin{array}{c}
\phi_{11}
\\
\phi_{21}
\\
\vdots
\\
\phi_{n1}
\end{array}\right)~,~
{\Phi}_{2}(x)=\left(
\begin{array}{c}
\phi_{12}
\\
\phi_{22}
\\
\vdots
\\
\phi_{n2}
\end{array}\right)~,\ldots,~
{\Phi}_{n}(x)=\left(
\begin{array}{c}
\phi_{1n}
\\
\phi_{2n}
\\
\vdots
\\
\phi_{nn}
\end{array}\right)~.
\end{equation}

The $n\times n$ determinant

\begin{equation}\label{wronskdef}
\left|
\begin{array}{cccc}
\phi_{11} &\phi_{12} &\ldots &\phi_{1n}
\\
\phi_{21} &\phi_{22} &\ldots &\phi_{2n}
\\
\vdots &\vdots & &\vdots
\\
\phi_{n1} &\phi_{n2} &\ldots &\phi_{nn}
\end{array}\right|
\end{equation}
is called the Wronskian of the $n$ vector functions ${\Phi}_{1}$,
${\Phi}_{2}$, $\ldots$ ,${\Phi}_{n}$ defined by (\ref{nvecdef}).
We will denote it by $W({\Phi}_{1}, {\Phi}_{2}, \ldots
,{\Phi}_{n})$ and its value at the point $x$ by $W({\Phi}_{1},
{\Phi}_{2}, \ldots ,{\Phi}_{n})(x)$.

Let the vector functions ${\Phi}_{1}$, ${\Phi}_{2}$, $\ldots$
,${\Phi}_{n}$ defined by (\ref{nvecdef}) be $n$ solutions of the
homogeneous linear vector differential equation (\ref{y0}) on the
real interval $[a,b]$. These $n$ solutions ${\Phi}_{1}$,
${\Phi}_{2}$, $\ldots$ ,${\Phi}_{n}$ of (\ref{y0}) are linearly
independent on $[a,b]$ if and only if
\begin{equation}\label{wronskneq0}
W({\Phi}_{1}, {\Phi}_{2}, \ldots ,{\Phi}_{n})(x)\neq0
\end{equation}
for all $x\in[a,b]$. A set of $n$ linearly independent solutions
of (\ref{y0}) is called a fundamental set of solutions of
(\ref{y0}). If the vector functions ${\Phi}_{1}$, ${\Phi}_{2}$,
$\ldots$ ,${\Phi}_{n}$ defined by (\ref{nvecdef}) make up a
fundamental set of solutions of (\ref{y0}), then the $n\times n$
square matrix
\begin{equation}\label{funmatrix}
{\bf M}(x)=\left(
\begin{array}{cccc}
\phi_{11}(x) &\phi_{12}(x) &\ldots &\phi_{1n}(x)
\\
\phi_{21}(x) &\phi_{22}(x) &\ldots &\phi_{2n}(x)
\\
\vdots &\vdots & &\vdots
\\
\phi_{n1}(x) &\phi_{n2}(x) &\ldots &\phi_{nn}(x)
\end{array}\right)
\end{equation}
is a fundamental matrix of (\ref{y0}). Furthermore, let ${\tilde
{\Phi}}_0(x)$ be an arbitrary solution of (\ref{y0}) on the real
interval $[a,b]$. Then there exists a suitable constant
vector\footnote{We denote the transpose of ${\bf A}$ by ${\bf
A}^{T}$, where {\bf A} being any matrix.}
\begin{equation}\label{convec}
{\bf c}=(c_1,c_2,\ldots,c_n)^T
\end{equation}
such that
\begin{equation}\label{solhom}
{\tilde {\Phi}}_0(x)={\bf M}(x){\bf c}
\end{equation}
on [a,b].

Inserting the solution (\ref{solhom}) of ${\bf y}_{0}$ into the
equation (\ref{y1}) yields
\begin{equation}
\label{y1eqn}  {\bf{y}}^{\prime}_{1}= {\bf{H}}(x){\bf{y}}_{1}+
{\bf{Z}}(x){\bf M}(x){\bf c}~.
\end{equation}
We shall now proceed to obtain a solution of above equation by
variation of parameters. It has been proven~\cite{ros80} that
${\tilde {\Phi}}_1$ defined by
\begin{equation}
\label{y1sol}  {\tilde {\Phi}}_1(x)= {\bf
M}(x)\int^{x}_{x_{0}}{\bf M}^{-1}(t) {\bf{Z}}(t){\bf M}(t){\bf
c}dt~,
\end{equation}
where $x_{0}\in[a,b]$, is a solution of the nonhomogeneous linear
vector differential equation (\ref{y1eqn}) on $[a,b]$.

We have shown that ${\tilde {\Phi}}_1$ is a solution of
differential equation (\ref{y1}). Directly inserting this solution
into the differential equation (\ref{y2}) yields
\begin{equation}
\label{y2eqn}  {\bf{y}}^{\prime}_{2}= {\bf{H}}(x){\bf{y}}_{2}+
{\bf{Z}}(x){\tilde {\Phi}}_1(x)~.
\end{equation}
We shall also proceed to obtain a solution of above equation by
variation of parameters. It has been proven that ${\tilde
{\Phi}}_2$ defined by
\begin{equation}
\label{y2sol}  {\tilde {\Phi}}_2(x)= {\bf
M}(x)\int^{x}_{x_{0}}{\bf M}^{-1}(t_2) {\bf{Z}}(t_2){\bf
M}(t_2)\int^{t_2}_{x_{0}}{\bf M}^{-1}(t_1) {\bf{Z}}(t_1){\bf
M}(t_1){\bf c}dt_1dt_2~,
\end{equation}
where $x_{0}\in[a,b]$, is a solution of the nonhomogeneous linear
vector differential equation (\ref{y2eqn}) on $[a,b]$.

According to above argument, we can acquire ${\tilde
{\Phi}}_3(x),{\tilde {\Phi}}_4(x),\ldots$ etc term by term. Now we
shall assume that the solutions ${\tilde {\Phi}}_0(x),{\tilde
{\Phi}}_1(x),\ldots,{\tilde {\Phi}}_{m-1}(x)$ have been acquired.
Then, from the recursion formula (\ref{ym}), we can simplify the
differential equation satisfied by ${\bf{y}}_{m}$
\begin{equation}
\label{ymeqn}  {\bf{y}}^{\prime}_{m}= {\bf{H}}(x){\bf{y}}_{m}+
{\bf{Z}}(x){\tilde {\Phi}}_{m-1}(x)~.
\end{equation}
It is easy to prove that ${\tilde {\Phi}}_m$ defined by
\begin{equation}
\label{ymsol} \begin{array}{l}\displaystyle
 {\tilde {\Phi}}_m(x)=
{\bf M}_x\int^{x}_{x_{0}}{\bf M}^{-1}_{t_m} {\bf{Z}}_{t_m}{\bf
M}_{t_m}\int^{t_m}_{x_{0}}{\bf M}^{-1}_{t_{m-1}}
\\\displaystyle~~~~~~~~~~~~~
\cdots{\bf M}_{t_2}\int^{t_2}_{x_{0}}{\bf M}^{-1}_{t_1}
{\bf{Z}}_{t_1}{\bf M}_{t_1}{\bf c}dt_1 dt_2 \cdots dt_m~,
\end{array}
\end{equation}
where ${\bf M}_t\equiv {\bf M}(t)$ and $x_{0}\in[a,b]$, is a
solution of the nonhomogeneous linear vector differential equation
(\ref{ym}) on $[a,b]$.

Hence the above reasoning demonstrate that ${\tilde {\Phi}}_j$
defined by
\begin{equation}
\label{yjsol} \begin{array}{l}\displaystyle
 {\tilde {\Phi}}_j(x)=
{\bf M}_x\int^{x}_{x_{0}}{\bf M}^{-1}_{t_j} {\bf{Z}}_{t_j}{\bf
M}_{t_j}\int^{t_j}_{x_{0}}{\bf M}^{-1}_{t_{j-1}}\cdots
\\\displaystyle~~~~~~~~~~~
{\bf M}_{t_2}\int^{t_2}_{x_{0}}{\bf M}^{-1}_{t_1}
{\bf{Z}}_{t_1}{\bf M}_{t_1}{\bf c}dt_1 dt_2 \cdots
dt_j~,~~~~j=1,2,3,\ldots
\end{array}
\end{equation}
where $x_{0}\in[a,b]$, is a solution of the nonhomogeneous linear
vector differential equation (\ref{con1}) on $[a,b]$.

It is well known that ${\Phi}_0(x)$ defined by
\begin{equation}
\label{y0sol}  {\Phi}_0(x)= {\bf M}(x)\int^{x}_{x_{0}}{\bf
M}^{-1}(t) {\bf F}(t)dt~,
\end{equation}
where $x_{0}\in[a,b]$, is a solution of the nonhomogeneous linear
vector differential equation
\begin{equation}
\label{nonhomy0} {\bf{y}}^{\prime}_{0}=
{\bf{H}}(x){\bf{y}}_{0}+{\bf{F}}(x)
\end{equation}
on the real interval $[a,b]$. A theorem~\cite{ros80} has been
proven that an arbitrary solution ${\bar {\Phi}}_0(x)$ of the
nonhomogeneous differential equation (\ref{nonhomy0}) is of the
form
\begin{equation}
\label{gensol} {\bar{\Phi}}_{0}(x)={\Phi}_0(x)+{\bf M}(x){\bf c}
\end{equation}
for a suitable choice of $c_{1},c_{2},\ldots,c_{n}$. According to
the same logic chain as what has been explained in the study of
the homogeneous linear vector differential equation
(\ref{dieqhom}), A solution of the nonhomogeneous linear vector
differential equation (\ref{con2}) on $[a,b]$ can easily be
acquired, that is
\begin{equation}
\label{yjsolnon} \begin{array}{l}\displaystyle
 {\bar {\Phi}}_j(x)=
{\bf M}_x\int^{x}_{x_{0}}{\bf M}^{-1}_{t_j} {\bf{Z}}_{t_j}{\bf
M}_{t_j}\int^{t_j}_{x_{0}}{\bf M}^{-1}_{t_{j-1}}\cdots
\\\displaystyle~~~~~~~
~~~~~ {\bf M}_{t_2}\int^{t_2}_{x_{0}}{\bf M}^{-1}_{t_1}
{\bf{Z}}_{t_1} \left( {\bf M}_{t_1}{\bf c}+{\bf
M}_{t_1}\int^{t_1}_{x_{0}}{\bf M}^{-1}_{t} {\bf{F}}_{t} dt \right)
dt_1 dt_2 \cdots dt_j~,
\\\displaystyle~~~~~~~
~~~~~j=1,2,3,\ldots
\end{array}
\end{equation}
where $x_{0}\in[a,b]$.

\section{The solution of the first-order
linear vector differential equation}

{\bf THEOREM}~~ the special series is the solution of the
first-order linear vector differential equation if the infinite
special series converges on the real interval $[a,b]$.

{\bf proof.} Consider the partial sums of vectors~\cite{vpr00}
\begin{equation}\nonumber
\begin{array}{l}
{\bf{S}}_{0}={\bf{y}}_{0}~,
\\
{\bf{S}}_{1}={\bf{y}}_{0}+{\bf{y}}_{1}~,
\\
{\bf{S}}_{2}={\bf{y}}_{0}+{\bf{y}}_{1}+{\bf{y}}_{2}~,
\\
\cdots~.
\end{array}
\end{equation}
Then ${\bf{S}}_{l}$, the $l$th partial sum of vectors, is given by
\begin{equation}\label{psn}
{\bf{S}}_{l}={\bf{y}}_{0}+{\bf{y}}_{1}+{\bf{y}}_{2}+\cdots+{\bf{y}}_{l}=\sum^{l}_{k=0}{\bf{y}}_{k}~.
\end{equation}
In the case of homogeneous linear vector differential equation
(\ref{dieqhom}), from the equations (\ref{con1}), we directly
acquire
\begin{equation}\label{psequ1}
{\bf{S}}^{\prime}_{l}= {\bf{H}}(x){\bf{S}}_{l}+
{\bf{Z}}(x){\bf{S}}_{l-1}~,~~~~l\geq 1~.
\end{equation}
In the case of nonhomogeneous linear vector differential equation
(\ref{dieq}), summing all the equations in (\ref{con2}) for the
set $0\leq j\leq l$, we obtain
\begin{equation}\label{psequ2}
{\bf{S}}^{\prime}_{l}= {\bf{H}}(x){\bf{S}}_{l}+
{\bf{Z}}(x){\bf{S}}_{l-1}+{\bf{F}}(x)~,~~~~l\geq 1~.
\end{equation}
Since the integer $l$ is unrelated with variable $x$, the limit
with respect to $l$ and the differential with respect to $x$ are
commutative
\begin{equation}\label{ldcom}
\lim_{l\to \infty} {\bf{S}}^{\prime}_{l} =(\lim_{l\to \infty}
{\bf{S}}_{l})^{\prime}~,
\end{equation}
With the help of equation (\ref{ldcom}), seeking for the limit of
the partial sum of vectors, for the homogeneous linear systems,
from the equation (\ref{psequ1}), we obtain
\begin{equation}\label{pslimt1}
(\lim_{l\to \infty} {\bf{S}}_{l})^{\prime}= {\bf{H}}(x)\lim_{l\to
\infty}{\bf{S}}_{l}+ {\bf{Z}}(x)\lim_{l\to \infty}{\bf{S}}_{l-1}~,
\end{equation}
and for the nonhomogeneous linear systems, seeking for the limit
of the partial sum of vectors in the equation (\ref{psequ2})
yields
\begin{equation}\label{pslimt2}
(\lim_{l\to \infty} {\bf{S}}_{l})^{\prime}= {\bf{H}}(x)\lim_{l\to
\infty}{\bf{S}}_{l}+ {\bf{Z}}(x)\lim_{l\to
\infty}{\bf{S}}_{l-1}+{\bf{F}}(x)~.
\end{equation}
Assuming that the infinite special series converges on the real
interval $[a,b]$, then
\begin{equation}\label{limit}
\lim_{l\to \infty} {\bf{S}}_{l} =\lim_{l\to \infty}
{\bf{S}}_{l-1}={\bf{S}}~,
\end{equation}
where ${\bf{S}}$ is a finite function on $[a,b]$. Therefore,
inserting the above expression into the equation (\ref{pslimt1})
yields
\begin{equation}
\label{sdif1}  {\bf{S}}^{\prime}= {\bf{H}}(x){\bf{S}}+
{\bf{Z}}(x){\bf{S}}={\bf{A}}(x){\bf{S}} ~,
\end{equation}
that is, the limit ${\bf{S}}$ defined by (\ref{con1}) and
(\ref{limit}) satisfies the homogeneous linear vector differential
equation (\ref{dieqhom}). Furthermore, inserting the formula
(\ref{limit}) into the equation (\ref{pslimt2}) yields
\begin{equation}
\label{sdif2}  {\bf{S}}^{\prime}= {\bf{H}}(x){\bf{S}}+
{\bf{Z}}(x){\bf{S}}+{\bf{F}}(x)={\bf{A}}(x){\bf{S}}+{\bf{F}}(x) ~.
\end{equation}
The above equation demonstrates that the limit ${\bf{S}}$ defined
by (\ref{con2}) and (\ref{limit}) satisfies the nonhomogeneous
linear vector differential equation (\ref{dieq}). Since the
infinite series converges and has sum ${\bf{S}}$ if the sequence
of partial sums converges to ${\bf{S}}$, the special series
converges to the solution of the first-order linear vector
differential equation.

Thus we have shown the validity of theorem.

{\bf CONCLUSION}~~If the special series converges, then:

(1) an arbitrary solution ${\tilde {\Phi}}(x)$ of the homogeneous
linear vector differential equation (\ref{dieqhom}) on $[a,b]$ can
be expressed as
\begin{equation}
\label{homsol}
 {\tilde {\Phi}}(x)=\sum_{j=0}^{\infty}{\tilde {\Phi}}_{j}(x)~,
\end{equation}
where ${\tilde {\Phi}}_{j}(x)$ is defined by (\ref{yjsol}); and

(2) an arbitrary solution ${\bar {\Phi}}(x)$ of the nonhomogeneous
linear vector differential equation (\ref{dieq}) on $[a,b]$ is of
the form
\begin{equation}
\label{nonhomsol}
 {\bar {\Phi}}(x)=\sum_{j=0}^{\infty}{\bar {\Phi}}_{j}(x)~,
\end{equation}
where ${\bar {\Phi}}_{j}(x)$ is defined by (\ref{yjsolnon}).

 {\bf Example}

Solve the system
\begin{equation}
\label{example}\left\{
\begin{array}{l}
y^{\prime}_{1}= (e^{-2x}-3e^{2x}+2) y_{1} + (e^{-2x}-9e^{2x}+3)
y_{2}
\\
y^{\prime}_{2}= (e^{2x}-e^{-2x}-1) y_{1} + (3e^{2x}-e^{-2x}-2)
y_{2}
\end{array}\right.
\end{equation}
on the real interval [0,1].

Clearly this is of the form (\ref{dieqhom2}), where
\begin{equation}
\label{coef} {\bf{H}}(x)=\left(
\begin{array}{cc}
2 & ~3
\\
-1 & ~-2
\end{array}\right),~~{\bf{Z}}(x)=\left(
\begin{array}{cc}
e^{-2x}-3e^{2x} & ~e^{-2x}-9e^{2x}
\\
e^{2x}-e^{-2x} & ~3e^{2x}-e^{-2x}
\end{array}\right)~.
\end{equation}
The corresponding homogeneous differential equation satisfied by
${\bf{y}}_{0}$ is
\begin{equation}
\label{basicde}\frac{d {\bf{y}}_{0}}{dx}=\left(
\begin{array}{cc}
2 & ~3
\\
-1 & ~-2
\end{array}\right){\bf{y}}_{0}~.
\end{equation}
we find that
\begin{equation}
\label{solv} \phi_{1}(x)=\left(
\begin{array}{c}
3e^{x}
\\
-e^{x}
\end{array}\right)~~~~{\rm and}~~~~\phi_{2}(x)=\left(
\begin{array}{c}
e^{-x}
\\
-e^{-x}
\end{array}\right)
\end{equation}
constitute a fundamental set (pair of linearly independent
solutions) of (\ref{basicde}). Thus a fundamental matrix of
(\ref{basicde}) is given by
\begin{equation}
\label{fmatrix} {\bf{M}}(x)=\left(
\begin{array}{cc}
3e^{x} & ~e^{-x}
\\
-e^{x} & ~-e^{-x}
\end{array}\right)~.
\end{equation}
From the fundamental matrix ${\bf{M}}(x)$, we find that
\begin{equation}
\label{finverse} {\bf{M}}^{-1}(x)=\frac{1}{2}\left(
\begin{array}{cc}
e^{-x} & ~e^{-x}
\\
-e^{x} & ~-3e^{x}
\end{array}\right)~.
\end{equation}
Thus the term ${\bf{M}}^{-1}(x){\bf{Z}}(x){\bf{M}}(x)$ in the
formula (\ref{yjsol}) becomes
\begin{equation}
\label{multipli} {\bf{M}}^{-1}(x){\bf{Z}}(x){\bf{M}}(x)=\left(
\begin{array}{cc}
0 & ~2
\\
2 & ~0
\end{array}\right)~.
\end{equation}
We have indicated that an arbitrary solution of (\ref{basicde})
can be presented as
\begin{equation}\label{exsolhom}
{\tilde {\Phi}}_0(x)={\bf M}(x)\left(
\begin{array}{c}
c_{1}
\\
c_{2}
\end{array}
\right)~,
\end{equation}
where $c_{1}$ and $c_{2}$ are constants. Inserting the equation
(\ref{multipli}) into the formula (\ref{y1sol}) yields\footnote{In
this example, we have selected $x_{0}=0$ for simplicity.}
\begin{equation}\label{exsolhom2}
{\tilde {\Phi}}_1(x)={\bf M}(x)\left(
\begin{array}{c}
2x c_{2}
\\
2x c_{1}
\end{array}
\right)~.
\end{equation}
One can acquire ${\tilde {\Phi}}_2(x)$, ${\tilde {\Phi}}_3(x)$,
$\ldots$ from the formula (\ref{yjsol}) term by term, that is
\begin{equation}\label{examsol}
{\tilde {\Phi}}_2(x)={\bf M}(x)\left(
\begin{array}{c}
\frac{(2x)^2}{2}c_{1}
\\
\frac{(2x)^2}{2}c_{2}
\end{array}
\right),~{\tilde {\Phi}}_3(x)={\bf M}(x)\left(
\begin{array}{c}
\frac{(2x)^3}{3!}c_{2}
\\
\frac{(2x)^3}{3!}c_{1}
\end{array}
\right),\ldots
\end{equation}
Obviously, the special series given by (\ref{exsolhom}),
(\ref{exsolhom2}), and (\ref{examsol}) is convergent. By theorem,
we know that
\begin{eqnarray}\nonumber
{\tilde {\Phi}}(x)&=&{\bf M}(x)\left(
\begin{array}{l}
c_{1}+2x
c_{2}+\frac{(2x)^2}{2!}c_{1}+\frac{(2x)^3}{3!}c_{2}+\frac{(2x)^4}{4!}c_{1}+\cdots
\\
c_{2}+2x
c_{1}+\frac{(2x)^2}{2!}c_{2}+\frac{(2x)^3}{3!}c_{1}+\frac{(2x)^4}{4!}c_{2}+\cdots
\end{array}
\right)
\\
\nonumber &=&\frac{1}{2}{\bf M}(x)\left(
\begin{array}{l}
e^{2x}(c_{1}+c_{2})+e^{-2x}(c_{1}-c_{2})
\\
e^{2x}(c_{2}+c_{1})+e^{-2x}(c_{2}-c_{1})
\end{array}
\right)
\end{eqnarray}
is a solution of the linear system (\ref{example}) for every real
number $x\in [0,1]$. Furthermore, we simplify the solution of
(\ref{example}) as follows
\begin{eqnarray}\nonumber
{\tilde {\Phi}}(x)=\left(
\begin{array}{c}
(-e^{-3x}+3e^{-x}+e^{x}+3e^{3x}){\tilde
c}_{1}+(e^{-3x}-3e^{-x}+e^{x}+3e^{3x}){\tilde c}_{2}
\\
(e^{-3x}-e^{-x}-e^{x}-e^{3x}){\tilde
c}_{1}+(-e^{-3x}+e^{-x}-e^{x}-e^{3x}){\tilde c}_{2}
\end{array}
\right)~,
\end{eqnarray}
where $2{\tilde c}_{1}=c_{1}$ and $2{\tilde c}_{2}=c_{2}$.
Inserting the above vector into the linear vector differential
equation (\ref{example}) really produces two identities.

An important fundamental property of a normal linear system
(\ref{dieq}) is its relationship to a single $n$th-order linear
differential equation in one unknown function. The so-called
normalized $n$th-order linear differential equation is of the form
\begin{equation}
\label{n-order}
\frac{d^{n}y}{dx^{n}}+a_{1}(x)\frac{d^{n-1}y}{dx^{n-1}}+\cdots+a_{n-1}(x)\frac{dy}{dx}+a_{n}(x)y=F(x)
~.
\end{equation}
It is well known that the single $n$th-order equation
(\ref{n-order}) can be transformed into a special case of the
normal linear system (\ref{dieq}) of $n$ equations in $n$ unknown
functions. Hence the method of special series can be used for
solving the normalized $n$th-order linear differential equation.

\section{Conclusions}

A special series has been introduced to yield the solution of the
linear vector differential equation. A group of recursion linear
vector differential equations satisfied by the terms of special
series has been  obtained and solved. We proved that the special
series is the solution of the first-order linear vector
differential equation if the infinite special series converges. We
have given an example to outline the procedure of our method. In
principle, the method of special series can be used for solving
any first-order linear vector differential equation if and only if
the considered special series converges.



\begin{thebibliography}{99}
\bibitem{hua05} X.~B.~Huang,
{\em ``Exact solutions of the Dirac equation in Robertson-Walker
space-time"}, gr-qc/0501077.

\bibitem{ros80} S.~L.~Ross,
{\em Introduction to ordinary differential equations}, third
edition, John Wiley $\&$ Sons, New York, 1980; S.~L.~Ross, {\em
Differential equations}, third edition, John Wiley $\&$ Sons, New
York, 1984.


\bibitem{wg00} Z.~X.~Wang, D.~R.~Guo, {\em Introduction to Special Function},
 Peking University Press, Beijing, 2000.

\bibitem{ww27} E.~T.~Whittaker, G.~N.~Watson,
{\em A course of modern analysis}, Cambridge, 1927.



\bibitem{vpr00} D.~Varberg, E.~J.~Purcell, S.~E.~Rigdon, {\em Calculus},
eighth edition, Prentice-Hall. Inc., New York, 2000.








\end{thebibliography}
\end{document}